
                             \def\flech{\rightarrow}

   \NoBlackBoxes \documentstyle{amsppt}

   \magnification=1100
    \loadbold
\topmatter
 \title An existence result of energy minimizer maps between riemannian
 polyhedra \endtitle

\leftheadtext{Taoufik Bouziane}

\rightheadtext{An existence result of energy minimizer maps}

\redefine\headmark#1{}
                          \author Taoufik Bouziane\endauthor

     \address  Taoufik BOUZIANE, The international centre for theoretical physics mathematic section,
      strada costiera 11, 34014 Trieste Italy  \endaddress

\email tbouzian\@ictp.trieste.it, btaoufik73\@hotmail.com\endemail

\abstract
In this paper, we prove the existence of energy
minimizers in each free homotopy class of maps between polyhedra
with target space without focal points. Our proof involves a
careful study of some geometric properties of riemannian polhyedra
without focal points. Among other things, we show that  on the
relevant polyhedra, there exists a convex supporting function.
 \endabstract

\keywords
 Geodesic space, Riemannian Polyhedra,Harmonic map, Focal point.
 E-mails: tbouzian\@ictp.trieste.it,
 btaoufik73\@hotmail.com
\endkeywords

\date  \enddate



   \endtopmatter


               \document

     \head{ 0. Introduction.}\endhead

    In the past decades, there has been a wide range of activity in the study of the existence
    of energy minimizers in various homotopy classes of maps between smooth riemannian
    manifolds, See [11], [15], [16], [18], [19], [20] and the references therein.
     In [11], Eells and Sampson obtained a fundamental theorem on the existence of harmonic
     maps in each free homotopy class of maps with target manifolds of non-positive sectional
     curvature which was generalized to the case of target manifolds without focal
     points by Xin [23]. In sharp contrast to the smooth
     case, very little is known for the case of singular spaces.
     In [17], Korevaar and Schoen expanded the theory of
     harmonic maps between smooth riemannian manifolds  to the
     case of maps  between certain singular spaces. For instance,
     admissible riemannian polyhedra are
     prototypes of the relevant singular spaces because these
     being both geodesic, harmonic (in the sense of Brelot cf. [10]
     ch 2), Dirichlet spaces and provide a wealth of examples as
     well. Let us mention here some examples of riemannian polyhedra (cf. [10] ch 8):

     \roster
     \item
     -Smooth riemannian manifolds (possibly with boundary).
\item
     -Triangulable riemannian lipshitz manifolds.
\item
     -Riemannian join of riemannian manifolds.
\item
     -Riemannian orbit spaces.
     \item
     -Riemannian orbifolds.
     \item
     -Conical singular riemannian spaces.
     \item
     -Normal analytic spaces with singularities.
     \item
     -Stratified riemannian spaces (or Thom spaces) satisfying
     Whitney's regularity condition. Etc...
     \endroster

     Our goal in the present paper is to show the existence of an
     energy minimizer in each free homotopy class of maps between
     riemannian polyhedra with target spaces without focal points
     in the sense of [4]. This result generalizes both Xin's result [23] and the new
     version of the Eells-Sampson's existence theorem [11] due
     to Eells and Fuglede (cf.[10] ch 11) where they obtained an existence
     theorem in the case of target spaces polyhedra of non-positive
     curvature in the sense of Alexandrov [6]. This
     generalization seems natural but it hides several
     difficulties which have to be solved by different approaches.
     One of these difficulties arises from the fact that in our case the absence
      of smoothness makes the Xin's methods [23] non-valid. Another difficulty, due to the fact
      that a riemannian polyhedron without focal points is not necessary of non-positive
      curvature, leads us making things differently from  Eells and Fuglede
      [11]. For instance, take for example a riemannian join of smooth riemannian manifold of positive sectional
      curvature without focal points and a riemannian manifold of non-positive
      curvature. This example has sense because Gulliver [14]
      has shown that there are manifolds without focal points of
      both signs of sectional curvature. In addition, Gulliver's result [14] implies that the
      nonexistence of focal points is weaker than non-positivity
      of the curvature.

     In order to state and prove our results we will give in section 1 some
     general preliminaries on geodesic spaces, riemannian
     polyhedra  and the energy of a map between riemannian
     polyhedra. In section 2, we will bring out some
     geometric properties of geodesic spaces which are due to the
     absence of focal points, we will also investigate the
     case of riemannian polyhedra without focal points in depth, and we will show the existence of
     a convex supporting function. The existence of such a function is the principal difference
     between our case and the Eells-Fuglede's case [10]. In their case,
     the square of the distance function to a geodesic is obviously
     convex supporting (consequence of the
     definition the non-positivity of the curvature). In our case the proof of such
     a fact is quite difficult. The geometric properties obtained in section 2,
     which are of self interest, are a subject of current
     investigations by the author. The results of such investigations will appear elsewhere.
     They are related to the dynamic of the  generalized geodesic flow in
     singular space. Lastly, section 3 is devoted to the existence of
     minimizing maps in free homotopy classes of maps between
     polyhedra.

                 \head { 1. preliminaries. }\endhead

This section is devoted to some preliminaries needed in the next
sections.

\subhead { 1.1. Geodesic spaces [2] [6] [7] [12] [13] }\endsubhead

Let $X$ be a metric space with metric $d$. A curve $c:I\flech X$
is called a {\it geodesic} if there is $v\geq 0$, called the
speed, such that every $t\in I$ has neighborhood $U\subset I$ with
$d(c(t_1),c(t_2))=v|t_1-t_2|$ for all $t_1,t_2\in U$. If the above
equality holds for all $t_1,t_2\in I$, then $c$ is called {\it
minimal geodesic}.

The space $X$ is called a {\it geodesic space} if every two points
in $X$ are connected by minimal geodesic. We assume from now on
that $X$ is a complete geodesic space.

A triangle $\Delta$ in $X$ is a triple $(\sigma_1, \sigma_2,
\sigma_3)$ of geodesic segments whose end points match in the
usual way. Denote by $H_k$ the simply connected complete surface
of constant Gauss curvature $k$. A {\it comparison triangle}
$\bar{\Delta}$ for a triangle $\Delta \subset X$ is a triangle in
$H_k$ with the same lengths of the sides as $\Delta$. A comparison
triangle in $H_k$ exists  and is unique up to congruence if the
lengths of sides of $\Delta$ satisfy the triangle inequality and,
in the case $k>0$, if the perimeter of $\Delta$ is
$<\frac{2\pi}{\sqrt k }$. Let $\bar{\Delta}= (\bar\sigma_1,
\bar\sigma_2, \bar\sigma_3)$ be a comparison triangle for
$\Delta=(\sigma_1, \sigma_2, \sigma_3)$, then for every point
$x\in \sigma_i$, $i=1,2,3$, we denote by $\bar x$ the unique point
on $\bar \sigma_i$ which lies at the same distances to the ends as
$x$.

Let $d$ denote the distance functions in both $X$ and $H_k$. A
triangle $\Delta$ in $X$ is $CAT_k$ {\it triangle} if the sides
satisfy the triangle inequality, the perimeter of $\Delta$ is
$<\frac{2\pi}{\sqrt k }$ for $k>0$, and if $d(x,y)\le d(\bar
x,\bar y)$, for every two points $x,y\in X$.

We say that $X$ has curvature at most $k$ and write $k_X\le k$ if
every point $x\in X$ has a neighborhood $U$ such that any triangle
in $X$ with vertices in $U$ and minimizing sides is $CAT_k$. Note
that we do not define $k_X$. If $X$ is a riemannian manifold, then
$k_X\le k$ iff $k$ is an upper bound for the sectional curvature
of $X$.

A geodesic space $X$ is called geodesically complete iff every
geodesic can be stretched in  two directions.

We say that a geodesic space $X$ is without conjugate points if
every two points in $X$ are connected by unique geodesic.

\subhead { 1.2. Orthogonality and focal point}\endsubhead

For more details on the study of focal points in geodesic space,
the reader can refer to [4].

\definition{Orthogonality}

$(X,d)$ will denote a complete geodesic space. Let $\sigma : \Bbb
R \flech X$ denote a geodesic and $\sigma _1:[a,b]\flech X$ a
minimal geodesic with a foot in $\sigma$ (i.e. $\sigma_1(a)\in
\sigma(\Bbb R)$).

The geodesic $\sigma_1$ is {\it orthogonal} to $\sigma$ if for all
$t\in [a,b]$, the point $\sigma_1(t)$ is locally of minimal
distance from $\sigma$.

In the case when for given geodesic $\sigma$ and a non-belonging
point $p$ there exists an orthogonal geodesic $\sigma'$ to
$\sigma$ and containing $p$, we will call the intersection point
between $\sigma$ and $\sigma'$ the {\it orthogonal projection
point} of $p$ on $\sigma$.

It is shown in the paper [4] that, on one hand, if the geodesic
$\sigma$ is minimal then  there always exists a realizing distance
orthogonal geodesic to $\sigma$ connecting every external point
$p$ (off $\sigma$)  to $\sigma$. On the other hand, if the space
$(X,d)$ is locally compact with non-null injectivity radius and
the geodesic $\sigma$ is minimal on every open interval with
length lower than the injectivity radius, then for every point $p$
off $\sigma$ and whose distance from $\sigma$ is not greater than
half of the injectivity radius, there exists a geodesic joining
orthogonally the point $p$ and the geodesic $\sigma$.

As corollaries, if the space $(X,d)$ is a simply connected $CAT_0$
space then for given geodesic $\sigma : \Bbb R \flech X$ and an
off point $p$ there always exists a realizing distance orthogonal
geodesic from $p$ to $\sigma$. When $X$ is $CAT_k$ for positive
constant $k$ then there always exists an orthogonal geodesic to
$\sigma$ from a point $p$ whose distance from $\sigma$ is not
greater than $\pi\over {2\sqrt k}$. In these last two cases the
angle between two orthogonal geodesics (in the sense of the
definition above) is always greater than or equal to $\pi\over 2$.

\enddefinition
 \definition{Focal points}

Let $(X,d)$ denote a complete geodesic space, $\sigma :\Bbb R
\flech X$ a geodesic and $p$ a point not belonging to the geodesic
$\sigma$.

The point $p$ is said to be a {\it focal point} of the geodesic
$\sigma$ or just a focal point of the space $X$, if there exists a
minimal geodesic variation $\tilde \sigma :]-\epsilon ,\epsilon [
\times[0,l]\flech X$ such that, if we note $\tilde \sigma
(t,s)=\sigma _t(s)$,  $\sigma _0$ is a minimal geodesic joining
$p$ to the point $q=\sigma (0)$ and for every $t\in ]-\epsilon
,\epsilon [$, $\sigma_t$ is a minimal geodesic containing
$\sigma(t)$, with the properties:

  \roster

  \item For every $t\in ]-\epsilon ,\epsilon[$,
        each geodesic $\sigma _t$ is orthogonal to $\sigma$.

        \item
       $\lim\limits _{t\to 0}{d(p,\sigma _t(l))\over {d(q,\sigma (t))}} =0$.
                             \endroster

 \enddefinition

 This definition was introduced in the article [4], us a
 natural generalization of the same notion in the smooth case.
 It is shown in the same paper that the Hadamard spaces are
 without a focal point.

\subhead { 1.3. Riemannian polyhedra}\endsubhead

\definition{ Riemannian admissible complexes ([3] [5] [6] [9] [22])}

 Let $K$ be a locally finite simplicial complex, endowed with a
 piecewise smooth riemannian metric $g$; i.e. $g$ is a family of
 smooth riemannian metrics $g_\Delta$ on simplices $\Delta$ of $K$
 such that the restriction $g_\Delta|\Delta'=g_{\Delta'}$ for any
 simplices $\Delta'$ and $\Delta$ with $\Delta'\subset \Delta$.

 Let $K$ be a finite dimensional simplicial complex which is connected
locally finite. A map $f$ from $[a,b]$ to $K$ is called a broken
geodesic if there is a subdivision $a=t_0<t_1<...<t_{p+1}=b$ such
that $f([t_i,t_{i+1}])$ is contained in some cell and the
restriction of f to $[t_i,t_{i+1}]$ is a geodesic inside that
cell. Then define the length of the broken geodesic map $f$ to be:
 $$ L(f)=\sum_{i=0}^{i=p} d(f(t_i),f(t_{i+1})).$$
 The length inside each cell being measured with respect to its metric.

 Then define $\tilde d(x,y)$, for every two points $x,y$ in $K$,
 to be the lower bound of the lengths of broken geodesics from $x$
 to $y$. $\tilde d$ is a pseudo-distance.

 If $K$ is connected and locally finite, then
$(K,\tilde d)$ is a length space and hence a geodesic space if is
complete (see also [5]).

A $l$-simplex in $K$ is called a {\it boundary simplex} if it is
adjacent to exactly one $l+1$ simplex. The complex $K$ is called
{\it boundaryless} if there are no boundary simplices in $K$.

The (open) {\it star} of an open simplex $\Delta^o$ (i.e. the
topological interior of $\Delta$ or the points of $\Delta$ not
belonging to any sub-face of $\Delta$, so if $\Delta$ is point
then $\Delta^o=\Delta$) of $K$ is defined as:
$$\text{$st(\Delta^o)=\bigcup \{ \Delta_i^o : \Delta_i$ is simplex
of $K$ with $\Delta_i\supset\Delta \}$ .}$$ The star $st(p)$ of
point $p$ is defined as the star of its {\it carrier}, the unique
open simplex $\Delta^o$ containing $p$. Every star is path
connected and contains the star of its points. In particular $K$
is locally path connected. The closure of any star is sub-complex.

\hskip 0 cm

 We say that the complex $K$ is {\it admissible}, if it is dimensionally homogeneous
and for every connected open subset $U$ of $K$, the open set
$U\setminus \{ U\cap \{\text{the $(k-2)$}-\text{skeleton}\} \}$ is
connected ($k$ is the dimension of $K$)(i.e. K is
$(n-1)$-chainable).

\hskip 0 cm

Let $x\in K$ a vertex of $K$ so that $x$ is in the $l$-simplex
$\Delta_{l}$. We view $\Delta_{l}$ as an affine simplex in $\Bbb
R^l$, that is $\Delta _l =\bigcap_{i=0}^l H_i$, where
$H_0,H_1,...,H_l$ are closed half spaces in general position, and
we suppose that $x$ is in the topological interior of $H_0$. The
riemannian metric $g_{\Delta_l}$ is the restriction to $\Delta_l$
of a smooth riemannian metric defined in an open neighborhood $V$
of $\Delta_l$ in $\Bbb R^l$. The intersection
$T_x\Delta_l=\bigcap_{i=1}^l H_i \subset T_xV$ is a cone with apex
$0\in T_xV$, and $g_{\Delta_l}(x)$ turns it into an euclidean
cone. Let $\Delta_m\subset \Delta_l$ ($m<l$) be another simplex
adjacent to $x$. Then, the face of $T_x\Delta_l$ corresponding to
$\Delta_m$ is isomorphic to $T_x\Delta_m$ and we view
$T_x\Delta_m$ as a subset of $T_x\Delta_l$.

Set $T_xK =\bigcup_{\Delta_i\ni x} T_x\Delta_i$, we call it the
{\it tangent cone} of $K$ at $x$. Let $S_x\Delta_l$ denote the
subset of all unit vectors in $T_x\Delta_l$ and set $S_x=S_xK
=\bigcup_{\Delta_i\ni x} S_x\Delta_i$. The set $S_x$ is called the
{\it link} of $x$ in $K$. If $\Delta_l$ is a simplex adjacent to
$x$, then $g_{\Delta_l}(x)$ defines a riemannian metric on the
$(l-1)$-simplex $S_x\Delta_l$. The family $g_x$ of riemannian
metrics $g_{\Delta_l}(x)$ turns $S_x\Delta_l$ into a simplicial
complex with a piecewise smooth riemannian metric such that the
simplices are spherical.

We call an admissible  connected locally finite simplicial
complex, endowed with a piecewise smooth riemannian metric, an
{\it admissible riemannian complex}.

 \enddefinition

\definition{ Riemannian polyhedron [10]}

We mean by {\it polyhedron} a connected locally compact separable
Hausdorff space $X$ for which there exists a simplicial complex
$K$ and homeomorphism $\theta : K \flech X$. Any such pair
$(K,\theta )$ is called a {\it triangulation} of $X$. The complex
$K$ is necessarily countable and locally finite (cf. [21] page
120) and the space $X$ is path connected and locally contractible.
The {\it dimension} of $X$ is by definition the dimension of $K$
and it is independent of the triangulation.

A {\it sub-polyhedron} of a polyhedron $X$ with given
triangulation $(K,\theta )$, is polyhedron $X'\subset X$ having as
a triangulation $(K',\theta|_{K'})$ where $K'$ is a subcomplex of
$K$ (i.e. $K'$ is complex whose vertices and simplexes are some of
those of $K$).

If $X$ is a polyhedron with specified triangulation $(K,\theta)$,
we shall speak of vertices, simplexes, $i-$skeletons or stars of
$X$ respectively of a space of links or tangent cones of $X$ as
the image under $\theta$ of vertices, simplexes, $i-$skeletons or
stars of $K$ respectively the image of space of links or tangent
cones of $K$. Thus our simplexes become compact subsets of $X$ and
the $i-$skeletons and stars become sub-polyhedrons of $X$.

If for given triangulation $(K,\theta)$ of the polyhedron $X$, the
homeomorphism $\theta$ is locally bi-lipschitz then $X$ is said
{\it Lip polyhedron} and $\theta$ {\it Lip homeomorphism}.

A {\it null set} in a Lip polyhedron $X$ is a set $Z\subset X$
such that $Z$ meets every maximal simplex $\Delta$, relative to a
triangulation $(K,\theta)$ (hence any) in set whose pre-image
under $\theta$ has $n-$dimensional Lebesgue measure $0$,
$n=dim\Delta$. Note that 'almost everywhere' (a.e.) means
everywhere except in some null set.

A {\it Riemannian polyhedron} $X=(X,g)$  is defined as a Lip
polyhedron $X$ with a specified triangulation $(K,\theta)$ such
that K is a simplicial complex endowed with a covariant bounded
measurable riemannian metric tensor $g$, satisfying the
ellipticity condition below. In fact, suppose that $X$ has
homogeneous dimension $n$ and choose a measurable riemannian
metric $g_\Delta$ on the open euclidean $n-$simplex
$\theta^{-1}(\Delta^o)$ of $K$. In terms of euclidean coordinates
$\{x_1,...,x_n\}$ of points $x=\theta^{-1}(p)$, $g_\Delta$ thus
assigns to almost every point $p\in \Delta^o$ (or $x$), an
$n\times n$ symmetric positive definite matrix $g_\Delta =
(g_{ij}^\Delta(x))_{i,j=1,...,n}$ with measurable real entries and
there is a constant $\Lambda_\Delta
>0$ such that (ellipticity condition):
$$\Lambda_\Delta^{-2}\sum_{i=0}^{i=n}(\xi^i)^2\le \sum_{i,j}
g^\Delta_{ij}(x)\xi^i\xi^j\le\Lambda_\Delta^2\sum_{i=0}^{i=n}(\xi^i)^2$$
for $a.e.$ $x\in\theta^{-1}(\Delta^o)$ and every
$\xi=(\xi^1,...,\xi^n) \in \Bbb R^n$. This condition amounts to
the components of $g_\Delta$ being bounded and it is independent
not only of the choice of the euclidean frame on
$\theta^{-1}(\Delta^o)$ but also of the chosen triangulation.

For simplicity of statements we shall sometimes require that,
relative to a fixed triangulation $(K,\theta)$ of riemannian
polyhedron $X$ (uniform ellipticity condition),
$$\text{$\Lambda$ $:=$
sup$\{\Lambda_\Delta:\Delta$ is simplex of $X\}<\infty$ .}$$

A riemannian polyhedron $X$ is said to be admissible if for a
fixed triangulation $(K,\theta)$ (hence any) the riemannian
simplicial complex $K$ is admissible.

There is a natural question we can ask about riemannian polyhedra:
Is the theorem of Gromov-Nash still true in the case of riemannian
polyhedra? In general, if we don't put more conditions on the
polyhedron, the answer to the question is no. In fact a
non-differentiable triangulable riemannian Lipschitz manifold is
an admissible riemannian polyhedron and, De Cecco and Palmieri [8]
showed that certain of these polyhedra are not isometrically
embeddable in any euclidean space (and therefore not in any smooth
riemannian manifold). But we know that finite dimensional Lip
polyhedron is affinely Lip embedded in some finite dimensional
euclidean space.

We underline that (for simplicity) the given definition of a
riemannian polyhedron $(X,g)$ contains already the fact (because
of the definition above of the riemannian admissible complex) that
the metric $g$ is {\it continuous} relative to some (hence any)
triangulation (i.e. for every maximal simplex $\Delta$ the metric
$g_\Delta$ is continuous up to the boundary). This fact is
sometimes in the literature omitted. The polyhedron is said to be
simplexwise smooth if relative to some triangulation $(K,\theta)$
(and hence any), the complex $K$ is simplexwise smooth. Both
continuity and simplexwise smoothness are preserved under
subdivision.

In the case of a general bounded measurable riemannian metric $g$
on $X$, we often consider, in addition to $g$, the {\it euclidean
riemannian metric} $g^e$ on the Lip polyhedron $X$ with a
specified triangulation $(K,\theta)$. For each simplex $\Delta$,
$g^e_\Delta$ is defined in terms of euclidean frame on
$\theta^{-1}(\Delta^o)$ as above by unitmatrix $(\delta_{ij})$.
Thus $g^e$ is by no means covariantly defined and should be
regarded as a mere reference metric on the triangulated polyhedron
$X$.

Relative to a given triangulation $(K,\theta)$ of an
$n-$dimensional riemannian polyhedron $(X,g)$ (not necessarily
admissible), we have on $X$ the distance function $e$ induced by
the euclidean distance on the euclidean space $V$ in which $K$ is
affinely Lip embedded. This distance $e$ is not intrinsic but it
will play an auxiliary role in defining an equivalent distance
$d_X$ as follows:

Let $\frak Z$ denote the collection of all null sets of $X$. For
given triangulation $(K,\theta)$ consider the set $Z_K\subset
\frak Z $ obtained from $X$ by removing from each maximal simplex
$\Delta$ in $X$ those points of $\Delta^o$ which are Lebesgue
points for $g_\Delta$. For $x,y \in X$ and any $Z\in \frak Z$ such
that $Z\subset Z_K$ we set:
$$\text{$d_X(x,y)=\sup \Sb {Z\in \frak Z}\\{Z\supset Z_K}\endSb
 \inf \Sb {\gamma}\\{\gamma(a)=x, \gamma(b)=y}\endSb \{ L_K(\gamma)$: $\gamma$ is Lip
continuous path and transversal to $Z\}$,}$$ where $L_K(\gamma)$
is  the length of the path $\gamma$ defined as:
$$\text{$L_K(\gamma)= \sum \Sb {\Delta\subset X}\\{}\endSb
 \int_{\gamma^{-1}(\Delta^o)} \sqrt{(g_{ij}^\Delta \circ \theta^{-1} \circ \gamma)
 \dot{\gamma}^i \dot{\gamma}^j } $, the sum is over all simplexes meeting $\gamma$.}$$

 It is shown in [10] that the distance $d_X$ is
 intrinsic, in particular it is independent of the chosen triangulation
 and it is equivalent to the euclidean distance $e$ (due to the Lip
 affinely and homeomorphically embedding of $X$ in some euclidean space $V$).

\enddefinition

\subhead { 1.4. Energy of maps}\endsubhead

The concept of energy in the case of a map of riemannian domain
into an arbitrary metric space $Y$ was defined and investigated by
Korevaar and Shoen [17]. Later this concept was extended by Eells
and Fuglede [10] to the case of map from an admissible riemannian
polyhedron $X$ with simplexwise smooth riemannian metric. Thus,
the energy $E(\varphi)$ of a map $\varphi$ from $X$ to the space
$Y$ is defined as the limit of suitable approximate energy
expressed in terms of the distance function $d_Y$ of $Y$.

 It is shown in [10] that the maps $\varphi : X\flech Y$ of
 finite energy are precisely those quasicontinuous (i.e.
 has a continuous restriction to closed sets, whose complements
 have arbitrarily small capacity, (cf. [10] page 153) whose restriction to each top dimensional simplex of $X$
 has finite energy in the sense of Korevaar-Schoen, and
 $E(\varphi)$ is the sum of the energies of these restrictions.

 Now, let $(X,g)$ be an admissible $m-$dimensional
 riemannian polyhedron with simplexwise smooth riemannian metric.
 It is not required that $g$ is continuous across lower
 dimensional simplexes. The target $(Y,d_Y)$ is an arbitrary
 metric space.

 Denote $L^2_{loc}(X,Y)$ the space of all $\mu_g-$measurable
 ($\mu_g$ the volume measure of $g$)
 maps $\varphi :X\flech Y$ having separable essential range and
 for which the map $d_Y(\varphi (.),q)\in L^2_{loc}(X,\mu_g)$
 (i.e. locally $\mu_g-$squared integrable) for some point $q$
 (hence by triangle inequality for any point). For $\varphi,\psi \in
 L^2_{loc}(X,Y)$ define their distance $D(\varphi,\psi)$ by:
 $$D^2(\varphi,\psi)= \int_X d_Y^2(\varphi (x),\psi(y)) d\mu_g(x).$$
 Two maps $\varphi,\psi \in L^2_{loc}(X,Y)$ are said to be {\it
 equivalent} if $D(\varphi,\psi)=0$, i.e. $\varphi(x)=\psi(x)$
 $\mu_g-$a.e. If the space $X$ is compact then $D(\varphi,\psi)<\infty$
 and $D$ is a metric on $L^2_{loc}(X,Y)=L^2(X,Y)$ and complete if
 the space $Y$ is complete [17].

 The {\it approximate energy density} of the map $\varphi\in
 L^2_{loc}(X,Y)$ is defined for $\epsilon >0$ by:
 $$e_\epsilon(\varphi)(x)=\int_{B_X(x,\epsilon)}\frac{d_Y^2(\varphi(x),\varphi(x'))}
 {\epsilon^{m+2}}d\mu_g(x').$$
 The function $e_\epsilon(\varphi)\ge 0$ is locally
 $\mu_g-$integrable.

 The {\it energy} $E(\varphi)$ of a map $\varphi$ of class
 $L^2_{loc}(X,Y)$ is:
 $$E(\varphi)=\sup_{f\in
 C_c(X,[0,1])}(\limsup_{\epsilon\rightarrow 0}\int_X f
 e_\epsilon(\varphi) d\mu_g),$$ where $C_c(X,[0,1])$ denotes the
 space of continuous functions from $X$ to the interval $[0,1]$
 with compact support.

A map $\varphi: X\flech Y$ is said to be {\it locally of finite
energy}, and we write $\varphi \in W^{1,2}_{loc}(X,Y)$, if
$E(\varphi|U)<\infty$ for every relatively compact domain
$U\subset X$, or equivalently if $X$ can be covered by domains
$U\subset X$ such that $E(\varphi|U)<\infty$.

For example (cf. [10] lemma 4.4), every Lip continuous map
$\varphi : X \flech Y$ is of class $W^{1,2}_{loc}(X,Y)$. In the
case when $X$ is compact $W^{1,2}_{loc}(X,Y)$ is denoted
$W^{1,2}(X,Y)$ the space of all maps of finite energy.

We can show  (cf. [10] theorem 9.1) that a map $\varphi \in
L^2_{loc}(X)$ is locally of finite energy iff there is a function
$e(\varphi)\in L^1_{loc}(X)$, named {\it energy density} of
$\varphi$, such that (weak convergence):
$$\text{ $\lim_{\epsilon\rightarrow 0}\int_X f e_\epsilon (\varphi) d\mu_g
=\int_X f e(\varphi) d\mu_g$,  for each $f\in C_c(X)$.}$$

\head { 2. Geodesic spaces without focal points. }\endhead

The aim of this section is to bring out some geometric properties
of a geodesic space, which are due to the absence of the focal
points. In particular, we will investigate in depth the example of
the riemannian polyhedra.

\subhead{ 2.1. Complete geodesic space without focal
points}\endsubhead

In this paragraph $(X,d)$ is a complete geodesic space. Our first
result is the following.

\proclaim{Theorem 2.1}

Suppose that the geodesic space $X$ is without focal points and it
is simply connected. Then $(X,d)$ is without conjugate points, in
the sense that for every pair of points $(p,q)\in X\times X$ there
exists a unique minimal geodesic $\sigma_{pq}$ connecting the
point $p$ to $q$.

\endproclaim

For the purpose of proving Theorem 2.1 we begin with the following
lemma:

\proclaim{ Lemma 2.2}

Under the same hypothesis of Theorem 2.1, let $\sigma : I\subseteq
\Bbb R \flech X$ be a geodesic and $p$ a point not belonging to
$\sigma$. If there exists a point $q\in \sigma$ which is an
orthogonal (see the paragraph 1.2) projection point of $p$ on
$\sigma$, then the point $q$ is unique.

\endproclaim

\demo{Proof}

In fact, the conclusion of the lemma means that the function
$t\mapsto L(t)=d^2(p,\sigma(t))$ reaches its minimum at most one
time.

Arguing by contradiction, then suppose that there exist $t_1\neq
t_2 \in I$ such that $\sigma(t_1)\neq\sigma(t_2)$ and both these
two points are orthogonal projection points of $p$ on $\sigma$.
Thus $d(\sigma(t_1),\sigma(t_2))> 0$ (because $\sigma$ is a
geodesic).

Now, let $\sigma_i$, $i=1,2$, be a minimal geodesic connecting the
point $p$ to $\sigma(t_i)$. According to the definition of
orthogonality, $\sigma_i$, for $i=1,2$, is orthogonal to the
geodesic $\sigma$.

We have supposed that our space $X$ is without focal points. This
fact is expressed by the following:
$$\text{$\forall s_1,s_2\ge 0$, $\exists K_{s_1s_2}>0$, such that,
$d(\sigma_1(s_1),\sigma_2(s_2))\ge
K_{s_1s_2}d(\sigma(t_1),\sigma(t_2))>0$.}$$

But this last assertion contradicts the fact that the geodesic
$\sigma_1$ meets the geodesic $\sigma_2$ at the point $p$.

This completes the proof of our result.

\hfill $\square$

\enddemo

Now, we are ready to prove Theorem 2.1.

\demo{Proof of Theorem 2.1}

 As in the theorem, let $(X,d)$ be a simply connected complete geodesic space
 without focal points. Suppose that there exists (at least) two
 points $p,q \in X$ such that, there is at least two distinct minimal geodesics
 $\sigma_1$ and $\sigma_2$, connecting them. Suppose that both $\sigma_1$
  and $\sigma_2$ are parameterized by the arc-length. Accordingly,
  there exists $s\in ]0,d(p,q)[$ such that
  $\sigma_1(s)\neq\sigma_2(s)$.

  Let $\gamma$ be a minimal geodesic connecting $\sigma_1(s)$ to
  $\sigma_2(s)$. Thus, by hypothesis we have
  $d(p,\sigma_1(s))=d(p,\sigma_2(s))$. By Lemma 2.2,
   there exists $t_0\in ]0,d(\sigma_1(s),\sigma_2(s))[$,
  such that $\gamma(t_0)$ is the unique orthogonal
  projection point of $p$ on the geodesic $\gamma$. So we have
  for every $t\neq t_0$, $d(p,\gamma(t_0))<d(p,\gamma(t))$. It
  follows from this last inequality the following:
  $$  \cases d(q,\gamma(t_0)) > d(q,\sigma_1(s))=d(q,\gamma(0))\\
 d(q,\gamma(t_0)) >
 d(q,\sigma_2(s))=d(q,\gamma(d(\sigma_1(s),\sigma_2(s))))& .\endcases$$

 Just now, look at the restrictions
 $\gamma|_{[0,t_0]}$ and
 $\gamma|_{[t_0,d(\sigma_1(s),\sigma_2(s))]}$, both are minimal
 geodesics. Therefore, using Lemma 2.2, we obtain two distinct orthogonal
 projection points of $q$, one of them is on the geodesic $\gamma([0,t_0])$
 the other is on the geodesic  $\gamma([t_0,d(\sigma_1(s),\sigma_2(s))])$.
 This last conclusion is in contradiction with Lemma 2.2 because the concatenation
 of the two geodesics is exactly $\gamma$ which is also a minimal
 geodesic. Our theorem is thereby proved.

\hfill $\square$

\enddemo

\subhead{ 2.2. Complete Riemannian polyhedra without focal
points}\endsubhead

This paragraph is devoted to a deep investigation of the geometry
of riemannian polyhedra without focal points. For simplicity of
statements we shall require that, our riemannian polyhedra are
simplexwise smooth. But the results of this paragraph are also
valid with mostly the same proofs if the riemannian polyhedra are
just Lip. First we will begin with some definitions.

Let $(X,d_X,g)$ be a complete riemannian polyhedron endowed with
simplexwise riemannian metric $g$ and $(K,\theta )$ a fixed
triangulation.

Recall that for each point $p\in X$ (or $\theta(p)\in K$), there
are well defined notions, the tangent cone over $p$ denoted $T_pX$
and the link over $p$ noted $S_pX$ which generalizes respectively
the tangent space and the unit tangent space  if $X$ is smooth
manifold (see the paragraph 1.3.).

Just now, we will suggest a generalization of the concepts of the
normal bundle  and the unit normal bundle of a geodesic in some
riemannian manifold, in the case of riemannian polyhedra.

\proclaim{Definitions 2.3}

Let $\sigma : I\subseteq \Bbb R \flech X$ be a geodesic, $p$ a
point belonging to $\sigma$ and $v\in T_pX$ a tangent vector.

\roster

\item $v$ is said {\it orthogonal} to the geodesic $\sigma$ iff,
there exists a geodesic $\gamma$ issuing from $p$ tangent to $v$
and orthogonal (see the definition above) to $\sigma$.

\item We name the set of all orthogonal vectors to $\sigma$ at the
point $p$ (could be empty), the {\it normal cone} of $\sigma$ over
$p$ and we denote it $\perp_p\sigma$.

\item We name the set of all unitary orthogonal vectors $u\in
\perp_p\sigma$, the {\it normal link} of $\sigma$ over $p$.

\endroster

\endproclaim

\proclaim{Remark 2.4}

As an immediate consequence we can derive from Theorem 2.1 that in
the case of riemannian polyhedra without focal points it make
sense to talk about a generalized exponential map $E_p : S_pX
\flech X$ which is an homeomorphism between the link over each
point $p$ and the space of all minimal geodesics deriving from
this point (because $X$ is locally without conjugate points).

\endproclaim

Next we prove a crucial geometric lemma which will play an
important role in all the following.

\proclaim{ Lemma 2.5}

Let $(X,d_X,g)$ be a complete riemannian polyhedron without focal
points and let $\sigma$ be a geodesic of $X$. Then for every point
$p$ belonging to $\sigma$, the spherical distance in the link
$S_pX$ between the two directions corresponding to the ingoing and
the outgoing of $\sigma$ at $p$, is greater or equal than $\pi$.

\endproclaim

\demo{Proof}

$(X,d_X,g)$ denotes a complete riemannian polyhedron without focal
points and $\sigma \subset X$ a geodesic. Suppose that there
exists a point $p= \sigma(0)$ where the conclusion of the lemma is
not valid.
 Let $\bigcup_{i=1}^n\Delta_i$ be a locale triangulation (we omit in the notation
  the homeomorphism of the triangulation) which
 contains $\sigma(0)$ in its (topological) interior. So there
 exists $\epsilon >0$ such that $\sigma(]-\epsilon, \epsilon[)\subset
 \bigcup_{i=1}^n\Delta_i$.

 As a consequence of what we suppose on $p$, the point $p$ necessarily belongs
  to the boundary of some simplex of the triangulation.
 Otherwise, $p$ will be in the (topological) interior of some
 simplex, but every open simplex is endowed with smooth metric and
 in this case, the distance between the two directions defined by
 $\sigma$ is equal to $\pi$ which is in contradiction with the
 hypothesis on the point $p$.

 Now suppose, that $p\in\partial \Delta_1\cap \partial \Delta_2$
 and note $v_1$, $v_2$ the two unitary tangent vectors to $\sigma$
 which are pointing respectively inside $\Delta_1$ and $\Delta_2$
 (these vectors are completely determined by the fact that the space is locally
 without focal points). So we have supposed that $dist(v_1,v_2)<\pi$
 ( $dist(v_1,v_2)$ means the spherical distance between $v_1$ and
 $v_2$). Thus, it should exist $\epsilon_0>0$ very close to $0$
  with $\sigma(\epsilon_0)\in \Delta_2$ and a neighborhood $U$
   of $p$ satisfying the radial uniqueness property such that:
 $$U \cap E(\perp_{\sigma(\epsilon_0)}\sigma) \cap E(\perp^{v_1}_{p}\sigma)\neq \emptyset,$$
with $E$ denoting the generalized exponential map and
$\perp^{v_1}_{p}\sigma$ the set of vectors $v\in
\perp_{\sigma(p)}\sigma$ which are orthogonal to $v_1$. This
contradicts the fact that the space $X$ is without focal points.

In fact, if such $\epsilon_0$ didn't exist, we will have for each
$t>0$:
 $$U \cap E(\perp_{\sigma(t)}\sigma)\cap E(\perp^{v_1}_{p}\sigma)= \emptyset,$$
 and by continuity (because the space is locally without focal points) we
 will have:
  $$U \cap E(\perp^{v_2}_p  \sigma)\cap E(\perp^{v_1}_{p}\sigma)= \emptyset,$$
  or, $$\text {$U \cap E(\perp^{v_2}_p  \sigma)\cap E(\perp^{v_1}_{p}\sigma)\supset \gamma$,
   where $\gamma$ is a minimal geodesic segment. }$$

   In fact, the last intersection cannot be a discrete set of points
   because that implies that $U$ contains conjugate points which
   is in contradiction with the radial uniqueness property.

  Both of the two last intersections  lead to the fact that the distance
  $dist(\perp^{v_2}_p  \sigma , \perp^{v_1}_{p}\sigma) =
  dist(v_1,v_2)-\pi$ is greater or equal to $0$, which leads to
  $dist(v_1,v_2\ge \pi$. This ends the proof of the lemma.

\hfill $\square$
\enddemo

Now, we are ready to prove the following theorem.

\proclaim{Theorem 2.6}

Let $(X,g,d)$ denote a complete simply connected riemannian
polyhedron without focal points, $\sigma: I\subseteq \Bbb R \flech
X$ a geodesic and $p\in X$. Then the function $L: t\mapsto
d^2(p,\sigma(t))$ is continuous and it is convex.

\endproclaim

\demo{Proof}

Firstly, the continuity of the function $L$ is a consequence of
the fact that the geodesic space $(X,d)$ is without conjugate
points (see Theorem 2.1).

Secondly, following the same argument used by Alexander and Bishop
in [1], where they show that a simply connected complete locally
convex geodesic space is globally convex, it is sufficient to show
that every point $x\in X$ admits an open convex neighborhood
$U_x$. In other terms, we just have to show the following: for
every $x\in X$ there is an open neighborhood $U_x$ such that,
every geodesic $\sigma$ with end points in $U_x$ belongs to $U_x$
and the function $L: t\mapsto d^2(x,\sigma(t))$ is convex.

Let $p$ be a point of the polyhedron $(X,g,d)$ and $(K,\theta)$ be
a fixed triangulation of $X$. In the following we will omit the
homeomorphism of the triangulation in our notations and so we will
do any distinction between the simplexes of $X$ and the simplexes
of $K$.

At first, we remark that for every $p\in X$, every real $r>0$ and
every geodesic $\sigma$ with ends in the open ball $B(p,r)$, with
center $p$ and ray $r$, is entirely contained in the ball
$B(p,r)$. In fact, the geodesic space $X$ is without focal points
so by the lemma of the last section (2.1.) we have:
$$\text{For every $t$, $d^2(p,\sigma(t))\le \sup\{d^2(p,$ first end
of $\sigma)$,$d^2(p,$ second end of $\sigma)$\}  .}$$

Second, there are two cases to investigate, the first one is when
the point $p$ is in the topological interior of some maximal
simplex and the second one is when the point $p$ is vertex (to the
triangulation $(K,\theta)$).

Suppose that $p$ is in the interior of the maximal simplex
$\Delta$. Then there exists a positive real $r_p>0$ such that the
open ball $B(p,r_p)$ with center $p$ and ray $r_p$ is contained in
$\Delta$. Thanks to the riemannian metric $g_{\Delta}$, the open
ball $B(p,r_p)$ can be thought of as sub-manifold of some smooth
riemannian manifold endowed with the riemannian metric
$g_{\Delta}$. Take now a geodesic $\sigma$ with end points in
$B(p,r_p)$ then it is  contained in the ball $B(p,r_p)$.

  The polyhedron $X$ is without focal points so the neighborhood (sub-manifold)
$B(p,r_p)$ is without focal points too. Thus, by a result of Xin
[23], the function $L$ is convex for every geodesic $\sigma$
contained in $B(p,r_p)$.

Now, look at the case where $p$ is vertex of $X$. Let $r_p$ be a
positive real such that the open ball $B(p,r_p)$ is included in
the open star $st(p)$ of $p$. Let $\sigma : [a,b] \flech X$ be a
geodesic of $B(p,r_p)$ and let $\bigcup_i \Delta_i^o$ (finite
union) denote the star of $p$. We know that there is a subdivision
$t_0=a,t_1,...,t_n=b$ such that each restriction
$\sigma|_{[t_i,t_{i+1}]}$ is a geodesic in the sense of smooth
riemannian geometry. So thanks to the result of Xin [23], the
question about the convexity of the function
$L(t)=d^2(p,\sigma(t))$ is asked when $\sigma$ transits from a
simplex $\Delta_i$ to a simplex $\Delta_{i+1}$ i.e. at the points
$t_i$.

Suppose that for fixed $t_j$ the function $L$ is not convex at
$t_j$. This hypothesis implies on one hand that $t_j$ is not the
minimum of the function $L$, because there is $\epsilon>0$ such
that both the restrictions $\sigma|_{[t_j-\epsilon,t_j]}$ and
$\sigma|_{[t_j,t_j +\epsilon]}$ are convex. So by taking an
$\epsilon$ smaller, the trace
$\sigma|_{[t_j-\epsilon,t_j+\epsilon]}$ is strictly monotone and
suppose it increasing. On the other hand, the non-convexity of $L$
at $t_j$ implies that the left derivative of $L$ at $t_j$ is
strictly greater than its right derivative. Let us now traduce
this last fact in terms of angles ; so let $\tau :[0,1]\times \Bbb
R \flech X$ be a map such that, for every $s\in \Bbb R$,
$\tau(.,s)=\tau_s(.)$ is the unique geodesic relating $p$ to
$\sigma(s)$ (because $X$ is globally without conjugate points).
Note $\theta_s^-$ and $\theta_s^+$ the left angle (to $s$)
respectively the right angle (to $s$) between the two geodesics
$\sigma$ and $\tau_s$. In a nutshell, if we traduce the
non-convexity of $L$ in term of angles we will have
$\pi-\theta_{t_j}^->\theta_{t_j}^+$ otherwise
$\theta_{t_j}^-+\theta_{t_j}^+<\pi$ which is in contradiction with
the last lemma (because $\sigma$ is a geodesic). So $L$ is still
convex at $t_j$. This ends the proof.

\hfill $\square$

\enddemo

\head{3. The existence theorem} \endhead

This section is devoted to the existence of minimizing maps in the
free homotopy classes  of maps between polyhedra. Henceforth all
polyhedra considered are supposed simplexwise smooth.

\proclaim{ Theorem 3.1}

Let $X$ and $Y$ be compact riemannian polyhedra. Suppose that $X$
is admissible and $Y$ is complete and without focal points.

Then every homotopy class $[u]$ of each continuous map $u$ between
the polyhedra $X$ and $Y$ has an energy (see the definition above)
minimizer relative to $[u]$.

\endproclaim

To prove this theorem we will adapt an original proof due to Eells
and Fuglede that we can find in [10], where they prove an
equivalent theorem in the case when the target polyhedron is
supposed of nonpositive curvature (in the sense of Alexandrov).
 In fact we will adjust the first step of their proof
 (because there are two steps in the proof of Eells and Fuglede)
 to our case, the second step remains the
same. But for the sake of completeness we will give all the proof
and just before we do some remarks.

\proclaim{Remarks 3.2}

An immediate consequence of Theorem 2.1, is that the universal
covering of a complete geodesic space without focal points is
contractible (because it is simple connected and without focal
points, see [1]). \vskip .3cm

Let $X$ and $Y$ be two locally finite polyhedra. If $Y$ has
contractible universal covering space, then the homotopy classes
of maps $u: X\flech Y$ are in natural bijective correspondence
$u\mapsto u_*$ with the conjugacy classes of homomorphisms $u_* :
\pi_1(X)\flech \pi_1(Y)$ of their fundamental groups (cf. [21] ch
8.1, theorem 9). \vskip .3cm

Let $X$ be an admissible compact riemannian polyhedron and
$(Y,d_y)$ be a complete geodesic space. Let a sequence
$(u_i)_i\subset W^{1,2}(X,Y)$ such that:
$$E(u_i)+ \int_X d_Y^2(u_i(x),q)d\mu_g(x)\le c,$$
for some constant $c$ and some fixed point $p\in Y$. Then
$(u_i)_i$ has a subsequence which converges in $L^2(X,Y)$ to a map
$u\in W^{1,2}(X,Y)$ (cf. [10] ch 9 Lemma 9.2).

\endproclaim

\vskip .5cm

After these remarks, we are actually ready to prove Theorem 3.1.

\demo{ Proof of Theorem 3.1}

Let $X$ and $Y$ be two riemannian polyhedra such that $X$ is
admissible and $Y$ is without focal points.

In the second section, and through the proof of Theorem 2.6, we
have shown that a riemannian polyhedron without focal points is
locally convex. Thus, in the case of the compact polyhedron $Y$
there exists a uniform $r>0$ such that for every point $y\in Y$,
the ball $B_Y(y,2r)$ is convex (i.e. $\forall \sigma: [a,b]\flech
Y$ with end points in $B_Y(y,2r)$, is completely contained in
$B_Y(y,2r)$ and the distance function $t\mapsto d^2(y,\sigma(t))$
is convex). So every geodesic contained in a ball $B_Y(y,2r)$ is
uniquely determined up to its end points.

So now, for given triangulations of both $X$ and $Y$ which are
compact, every continuous map $u:X\flech Y$ can be approximated
uniformly (because $X$ and $Y$ are compact) by a simplicial map
$u^{S}$ which is Lipschitz and hence of finite energy (see 1.4).
In addition, if we assume that $d_Y(u(x)-u^S(x))<r$ for all $x\in
X$, the simplicial map $u^S$ becomes homotopic to $u$. In fact,
finitely many  balls $B_Y(y_i,2r)$, $i=1,...,l$ cover the
polyhedron $Y$, so for given $i$ and $x\in U_i:=
u^{-1}(B_Y(y_i,2r))$ there is a unique minimal geodesic
$\sigma_x:I\subset \Bbb R \flech B_Y(y_i,2r)$ joining $u(x)$ to
$u^S(x)$ within the ball $B_Y(y_i,2r)$. Thus, the map $\tau:
I\times U_i\flech B_Y(y_i,2r)$ such that $\tau(t,x)=\sigma_x(t)$,
is continuous because geodesic segments in the convex balls
$B_Y(y_i,2r)$ vary continuously with their endpoints [1].

In summary, we have shown that for every continuous map $u:X\flech
Y$, there is a representative element of the homotopy class $[u]$
which is of finite energy.

Now, take an energy minimizing sequence $(u_i)_i$ of continuous
maps of finite energy in a homotopy class $[u]$. By the last
remark (see above) and the fact that, $d_Y(u_i(x),q)$ is bounded
by the diameter of $Y$ (because $Y$ is compact) and the polyhedron
$X$ is of finite volume ($X$ is compact), there is a
 subsequence,  noted always $(u_i)_i$, which converges in
 $L^2(X,Y)$ to an element $u\in W^{1,2}(X,Y)$. Every element $u_i$
 from this convergent subsequence $(u_i)_i$ lifts to a continuous
 map $\tilde u_i \in W^{1,2}_{loc}(\tilde X, \tilde Y)$ where
 $\tilde X$ and $\tilde Y$ are respectively the universal cover of
 $X$ and the universal cover of $Y$. In addition, every such map $\tilde u_i$
 is equivariant with respect to their fundamental  groups $\pi_1(X)$ and
 $\pi_1(Y)$ in the sense that, if $ (u_i)_*:\pi_1(X)\flech \pi_1(Y)$
 denote the homomorphism induced by the map $u_i$, then:
 $$\text{ $\tilde u_i\circ \gamma =(u_i)_*(\gamma)\circ \tilde
 u_i$ \hskip 1cm for all $\gamma\in \pi_1(X)$ \hskip 1cm  $(\star_i)$.}$$
 To normalize these lifted map we fix a point $x_0\in X_0$ and
 choose an image point noted $\tilde u_i(\tilde x_0)\in \tilde Y$ from
 the inverse image of the point $u_i(x_0)$ by the covering map (to the universal covering)
 of $Y$. The second remark above insures that the class $[u]$ is
 identified with the conjugacy class of homomorphism $u_* :
\pi_1(X)\flech \pi_1(Y)$ so independently of the choice of $\tilde
u_i(\tilde x_0)$, in our case of the minimizing subsequence
$(u_i)_i$ we can choose for all $i$ the same representative of the
class of the homeomorphism and let this choice be $(\tilde u_i)_*
= (\tilde u_1)_* $. Thus, the pointwise limit $\tilde u =
\lim_{i\rightarrow \infty} \tilde u_i$ satisfies the following:
$$\text{ $\tilde u\circ \gamma =(u_1)_*(\gamma)\circ \tilde u$
\hskip 1cm for all $\gamma\in \pi_1(X)$ \hskip 1cm $(\star)$.}$$

Now, we will use some algebraic topological arguments. Remember
that the fundamental group $\pi_1(X)$ acts isometrically and
simplicially on $\tilde X$ thus there exists  a compact set
$\tilde F\subset \tilde X$ called a {\it fundamental domain} of
$\pi_1(X)$ whose boundary $\partial \tilde F$ has measure $0$ and
each point of $\tilde X$ is $\pi_1(X)$-equivalent either to
exactly one point of the interior of $\tilde F$ or to at least one
point of $\partial \tilde F$. The fact that $X$ is compact implies
that the compact $\tilde F$ can be obtained as a suitable union of
maximal simplexes of $\tilde X$. Furthermore, $\tilde F$ is
contained in the interior $\tilde U$ of suitable union of
finitely, say $N$, many $\pi_1(X)$-translates of $\tilde F$.

Let $\frak E$ denote the class of all maps in $
W^{1,2}_{loc}(\tilde X, \tilde Y)$ which are equivariant as in
$(\star)$. So as we saw, the limit map $\tilde u =
\lim_{i\rightarrow \infty} \tilde u_i$ belongs to $\frak E$.

Actually, modify the above construction. Let $(\tilde u_i)_i$
denote a minimizing sequence for $\int_{\tilde F} e(\tilde u)$ in
the class $\frak E$, $e(\tilde u)$ denoting the energy density of
the map $\tilde u \in W^{1,2}_{loc}(\tilde X, \tilde Y)$. Thanks
to the equivariance equality $(\star)$, the sequence $(\tilde
u_i)_i$ is likewise minimizing for $E(\tilde u_{|\tilde U})= N
\int_{\tilde F} e(\tilde u)$. Now, by the third remark of 3.2.
applied to compact subsets of $\tilde U$, thus the
 sequence of traces $(\tilde u_{i_{|\tilde U}})$ converges in
 $L^2_{loc}(\tilde U, \tilde Y)$ and pointwise $a.e.$ in $\tilde
 U$ to some map $\tilde u_{\tilde U}\in W^{1,2}(\tilde U,\tilde Y)$
 which minimizes the energy of restrictions to $\tilde U$ of all maps
 belonging to the class $\frak E$. Consequently, and by $(\star)$
 again the sequence $(\tilde u_i)_i$ converges pointwise $a.e.$ in
 $\tilde X$ to an extension $\tilde u$ of the map $\tilde u_{\tilde
 U}$. Of course this new limit $\tilde u$ satisfies $(\star)$ but
 likewise it minimizes the integral
 $\int_{\tilde F} e(\tilde u) = N^{-1}E(\tilde u_{|\tilde U})$
 among all the restrictions to $\tilde F$ of maps of class $\frak
 E$. Furthermore, such minimizer  is also locally $E$-minimizing
 on $\tilde X$ ; indeed, every point  of $\tilde X$ has a
 relatively compact neighborhood $\tilde V$ such that $\tilde V
 \cap \gamma(\tilde V)= \emptyset$ for all $\gamma\in \pi_1(X)\setminus
 \{id\}$. So if an element $\tilde v\in W^{1,2}_{loc}(\tilde X, \tilde Y)$
 satisfies $\tilde v=\tilde u$ in $\tilde X\setminus \tilde V$,
 then the map $\tilde v^* :\tilde X\flech \tilde Y$ defined by
 $\tilde v^*(\gamma\tilde x)=\tilde v(\tilde x)$ for every $\tilde
 x\in \tilde V$ and $\gamma \in \pi_1(X)$, while $\tilde v^*
 =\tilde v$ elsewhere, belongs to the class $\frak E$, and
 satisfies $E(\tilde v_{|\tilde V})\ge E((\tilde v^*)_{|\tilde V})
 \ge E(\tilde u_{|\tilde  V})$.

 Furthermore the map $u:X\flech Y$ covered by  $\tilde u$ is in
 the class $[u]$ and minimizes the energy in $[u]$ ; indeed, any
 map $v\in [u]$ lifts to a map $\tilde v$ belonging to the class
 $\frak E$, and
 $E(v)=\int_{\tilde F} e(\tilde v)\ge \int_{\tilde F} e(\tilde
 u)=E(u)$ because $\tilde u$ is minimizing relative to the class
 $\frak E$. This ends the proof of the theorem.

\hfill $\square$

\enddemo

\subhead{Acknowledgement}\endsubhead

 The author would like to express his thanks to Prof. J
Eells for encouraging him to investigate this subject.

                                       \enddocument
\end